\documentclass[a4paper]{amsart}

\usepackage{amssymb}
\usepackage[arrow]{xy}
\usepackage{pb-diagram,pb-xy}
\usepackage{graphicx}

\theoremstyle{plain}
\newtheorem{theorem}{Theorem}[section]
\newtheorem{proposition}[theorem]{Proposition}

\newtheorem{lemma}[theorem]{Lemma}

\theoremstyle{definition}

\newtheorem{remark}[theorem]{Remark}

\def\Z{\mathbb{Z}}
\def\Q{\mathbb{Q}}
\def\R{\mathbb{R}}

\def\sR{\mathcal{R}}
\def\K{\mathcal{K}}
\def\U{\mathcal{U}}
\def\N{\mathcal{N}}

\def\Ker{\operatorname{Ker}}
\def\Im{\operatorname{Im}}
\def\Hom{\operatorname{Hom}}
\def\sign{\operatorname{sign}}

\def\tsg{g^{t}_*}
\def\ssg{g^{s}_*}

\def\to{\mathchoice{\longrightarrow}{\rightarrow}{\rightarrow}{\rightarrow}}
\makeatletter
\newcommand{\shortxra}[2][]{\ext@arrow 0359\rightarrowfill@{#1}{#2}}
\def\longrightarrowfill@{\arrowfill@\relbar\relbar\longrightarrow}
\newcommand{\longxra}[2][]{\ext@arrow 0359\longrightarrowfill@{#1}{#2}}
\renewcommand{\xrightarrow}[2][]{\mathchoice{\longxra[#1]{#2}}%
  {\shortxra[#1]{#2}}{\shortxra[#1]{#2}}{\shortxra[#1]{#2}}}
\makeatother

\makeatletter
\def\Nopagebreak{\@nobreaktrue\nopagebreak}
\makeatother

\begin{document}

\title{Topological minimal genus and $L^2$-signatures}

\author{Jae Choon Cha}

\address{Information and Communications University, Munji-dong,
  Yuseong-gu, Daejeon 305--732, Republic of Korea}

\email{jccha@icu.ac.kr}

\def\subjclassname{\textup{2000} Mathematics Subject Classification}
\expandafter\let\csname subjclassname@1991\endcsname=\subjclassname
\expandafter\let\csname subjclassname@2000\endcsname=\subjclassname


\begin{abstract}
  We obtain new lower bounds of the minimal genus of a locally flat
  surface representing a 2-dimensional homology class in a topological
  4-manifold with boundary, using the von Neumann-Cheeger-Gromov
  $\rho$-invariant.  As an application our results are employed to
  investigate the slice genus of knots.  We illustrate examples with
  arbitrarily large slice genus for which our lower bound is optimal
  but all previously known invariants vanish.
\end{abstract}

\maketitle

%

\section{Introduction and main results}

This paper concerns the problem of the minimal genus of a locally flat
embedded surface representing a given 2-dimensional homology class in
a topological 4-manifold.  Precisely, a locally flat closed surface
$\Sigma$ in a topological 4-manifold $W$ is said to represent
$\sigma\in H_2(W)$ if the fundamental class of $\Sigma$ is sent to
$\sigma$ under the map induced by the inclusion.  In this paper
manifolds are always oriented and surfaces are assumed to be
connected.

For topological 4-manifolds which are closed (or with boundaries
consisting of homology spheres), there are remarkable known results
which provide lower bounds of the minimal genus, including
Kervaire-Milnor \cite{Kervaire-Milnor:1961-1}, Hsiang-Szczarba
\cite{Hsiang-Szczarba:1971-1}, Rokhlin \cite{Rokhlin:1971-1}, and
Lee-Wilczy\'nski \cite{Lee-Wilczynski:1997-1,Lee-Wilczynski:2000-1}.
Basically these lower bounds are extracted by considering Rokhlin's
theorem and algebraic topology of finite cyclic branched coverings.
Also, interesting results on the smooth analogue of this problem have
been obtained using gauge theory for a certain class of smooth
4-manifolds.  Related works include results on the Thom conjecture and
adjunction inequality due to Kronheimer-Mrowka
\cite{Kronheimer-Mrowka:1993-1, Kronheimer-Mrowka:1995-1,
  Kronheimer-Mrowka:1995-2}, Morgan-Szab\'o-Taubes
\cite{Morgan-Szabo-Taubes:1996-1}, Kronheimer
\cite{Kronheimer:1999-1}, and Ozsv\'ath-Szab\'o
\cite{Ozsvath-Szabo:2000-1,Ozsvath-Szabo:2000-2}.  One may obtain a
lower bound in a 4-manifold with boundary when it embeds into another
4-manifold for which the above lower bound results can be applied
directly.  For example, the adjunction inequality is proved in Stein
4-manifolds via an embedding theorem due to Lisca-Mati\'c
\cite{Lisca-Matic:1997-1, Lisca-Matic:1998-1,
  Akbulut-Matveyev:1997-1}.

In this paper we focus on the minimal genus problem in a topological
4-manifold which has boundary with nontrivial homology.  Our results
give new lower bounds of the minimal genus, for homology classes from
the boundary, in terms of the von Neumann-Cheeger-Gromov
$\rho$-invariants of the boundary.  As an application we give lower
bounds of the slice genus of a knot.  Examples illustrate that the
$\rho$-invariants detect arbitrarily large minimal genus that all
previously known results do not.

\subsection*{Minimal second Betti number of a 4-dimensional bordism}

We obtain lower bounds of the minimal genus through the following
problem on 4-dimensional bordisms: what is the minimal second Betti
number of a topological null-bordism of a given closed 3-manifold
endowed with a group homomorphism of the fundamental group?  Our
principal result on this is as follows.  Let $\Gamma$ be a
poly-torsion-free-abelian (PTFA) group, i.e., $\Gamma$ admits a finite
length normal series $\{G_i\}$ with $G_i/G_{i+1}$ torsion-free
abelian.  It is known that there is a (skew-)field $\K$ of right
quotients of~$\Z\Gamma$.  Let $\sR$ be a subring of $\K$ which is a
PID containing~$\Z\Gamma$.  Then $\Gamma$ acts on the abelian group
$\K/\sR$ via right multiplication so that the semi-direct product
$(\K/\sR) \rtimes \Gamma$ is defined.

\begin{theorem}
  \label{theorem:lower-bound-of-minimal-second-betti-number}
  Suppose $M$ is a closed 3-manifold endowed with a homomorphism
  $\phi\colon \pi_1(M) \to \Gamma$, and $W$ is a topological
  4-manifold with boundary $M$ such that $\phi$ factors
  through~$\pi_1(W)$.  Then the followings hold:
  \begin{enumerate}
  \item
    The second Betti number $\beta_2(W)$ satisfies
    \[
    |\rho(M, \phi)| \le 2\beta_2(W)
    \]
    where $\rho(M, \phi)\in \R$ denotes the von Neumann-Cheeger-Gromov
    $\rho$-invariant of $M$ associated to~$\phi$.
  \item In addition, if the twisted homology $H_1(M;\sR)$ is
    $\sR$-torsion and not generated by any $\beta_2(W)$ elements, then
    there is a nontrivial submodule $P$ in $\Hom(H_1(M;\sR),\K/\sR)$
    such that every homomorphism in $P$ gives rise to a lift
    $\phi_1\colon \pi_1(M) \to (\K/\sR) \rtimes \Gamma $ of $\phi$
    which factors through~$\pi_1(W)$.
  \end{enumerate}
  
\end{theorem}

More detailed versions of Theorem
\ref{theorem:lower-bound-of-minimal-second-betti-number} (1) and (2)
are stated and proved in
Section~\ref{section:Betti-numbers-and-L2-signatures}
and~~\ref{section:extending-coefficient-systems}, respectively.  To
prove (1), we regard the $\rho$-invariant of $M$ as an $L^2$-signature
defect of $W$, and estimate the $L^2$-signature of $W$ in terms of the
$L^2$ and ordinary Betti number.  While (1) gives a lower bound of
$\beta_2(W)$ without using (2), further information may be obtained
when (2) is combined with~(1); note that (2) gives a sufficient
condition which implies that a certain ``bigger'' coefficient system
of $M$, namely~$\phi_1$, extends to~$W$.  In case that~$\phi_1$
extends, (1)~can be applied again to $\phi_1$ to obtain further lower
bounds of $\beta_2(W)$ (and possibly this process may be iterated).

This type of coefficient extension problem plays a crucial role in
earlier landmark works in knot theory, including Casson and
Gordon~\cite{Casson-Gordon:1986-1,Casson-Gordon:1978-1},
Gilmer~\cite{Gilmer:1982-1}, and in particular Cochran, Orr, and
Teichner~\cite{Cochran-Orr-Teichner:1999-1,Cochran-Orr-Teichner:2002-1},
from which
Theorem~\ref{theorem:lower-bound-of-minimal-second-betti-number} has
been directly motivated.  In
\cite{Cochran-Orr-Teichner:1999-1,Cochran-Orr-Teichner:2002-1} the
extension problem is investigated when $H_1(\partial W;\Q)\cong
H_1(W;\Q)\cong \Q$ and $W$ satisfies some geometric condition related
to the existence of a Whitney tower (such $W$ is called an
$(h)$-solution in~\cite{Cochran-Orr-Teichner:1999-1}).  In order to
deal with the extension problem without assuming these conditions, as
in
Theorem~\ref{theorem:lower-bound-of-minimal-second-betti-number}~(2),
we investigate the relationship of the Blanchfield linking form of $M$
and the intersection form of~$W$ over $\sR$-coefficients, and import
ideas from Gilmer's work~\cite{Gilmer:1982-1} on Casson-Gordon
invariants.

The following result relates the minimal second Betti number of
bordisms with a particular type of the minimal genus problem in a
4-manifold with boundary.  Suppose $W$ is a topological 4-manifold
with boundary $M$, $H_1(W)=0$, and $\sigma$ is a 2-dimensional
homology class contained in the image of $H_2(M) \to H_2(W)$.  In
Section~\ref{section:construction-of-bordisms-from-locally-flat-surfaces}
we will describe a homomorphism $\phi_\sigma\colon \pi_1(M) \to \Z$
determined by~$\sigma$.

\begin{proposition}
  \label{proposition:Z-bordism-from-embedded-surface}
  If $\phi_\sigma$ is nontrivial and there is a locally flat embedded
  surface of genus $g$ in $W$ representing $\sigma$, then there is a
  topological 4-manifold $V$ bounded by $M$ such that
  $\phi_\sigma\colon \pi_1(M) \to \Z$ factors through $\pi_1(V)$ and
  $\beta_2(V) = \beta_2(W)+2g-1$.
\end{proposition}

Consequently, lower bounds of $\beta_2(V)$ obtained by (possibly
repeatedly) applying
Theorem~\ref{theorem:lower-bound-of-minimal-second-betti-number} give
rise to lower bounds of the genus~$g$.

\subsection*{Slice genus of a knot}

As an application, we employ our results on the minimal genus problem
to investigate the slice genus of a knot $K$ in $S^3$.  The
\emph{topological slice genus} $\tsg(K)$ of $K$ is defined to be the
minimal genus of a locally flat surface $F$ properly embedded in $D^4$
in such a way that $\partial F=K$, viewing $S^3$ as the boundary
of~$D^4$.  The \emph{smooth slice genus} $\ssg(K)$ is defined
similarly, requiring $F$ to be a smooth submanifold of~$D^4$.
Obviously $\tsg(K) \le \ssg(K)$.

There are various known lower bounds of the slice genus.  Clearly any
obstruction to being a slice knot can be viewed as a lower bound of
the form (slice genus) $\ge 1$.  It is well known that some invariants
derived from a Seifert matrix, including the signature of a knot, can
be used to detect higher topological slice genus.  Gilmer showed that
Casson-Gordon invariants of a knot $K$ give further lower bounds
of~$\tsg(K)$~\cite{Gilmer:1982-1}.  For the smooth slice genus,
further results based on gauge theory are known.  For a special class
of knots which includes the torus knots, an optimal lower bound is
obtained as an application of the Thom conjecture due to
Mrowka-Kronheimer~\cite{Kronheimer-Mrowka:1993-1}.  For an arbitrarily
given knot $K$, the Thurston-Bennequin invariant (together with the
rotation invariant) of a Legendrian representation of $K$ is known to
give a lower bound of $\ssg(K)$, due to Rudolph~\cite{Rudolph:1995-1,
  Rudolph:1997-1}, Kronheimer-Mrowka,
Akbulut-Matveyev~\cite{Akbulut-Matveyev:1997-1}, and
Lisca-Mati\'c~\cite{Lisca-Matic:1997-1, Lisca-Matic:1998-1}.  More
recently, Ozsv\'ath-Szab\'o's
$\tau$-invariant~\cite{Ozsvath-Szabo:2003-1} and Rasmussen's
$s$-invariant~\cite{Rasmussen:2004-1} defined from knot homology
theories of Ozsv\'ath-Szab\'o and Khovanov have been known to give new
lower bounds of~$\ssg(K)$.

It is well-known that lower bounds of the slice genus can be obtained
through minimal genus problems in 4-manifolds with boundary; the slice
genus of a knot $K$ is bounded from below by the minimal genus for a
specific homology class in the 4-manifold obtained by attaching a
2-handle to the 4-ball along~$K$.  It follows that
Theorem~\ref{theorem:lower-bound-of-minimal-second-betti-number} and
Proposition~\ref{proposition:Z-bordism-from-embedded-surface} give
lower bounds of the slice genus in terms of the $\rho$-invariants.  In
fact, it turns out that this method gives us lower bounds of the genus
of a locally flat surface bounded by $K$ in a homology 4-ball with
boundary~$S^3$.  The following theorem illustrates that our lower
bounds from the $\rho$-invariants actually reveal new information; one
can detect arbitrarily large slice genus of knots that all the
previously known lower bounds fail to detect.

\begin{theorem}
  \label{theorem:main-example}
  For any positive integer $g$, there are infinitely many knots $K$
  with the following properties:
  \begin{enumerate}
  \item $\tsg(K)=\ssg(K)=g$.
  \item $K$ has a Seifert matrix of a slice knot.
  \item $K$ has vanishing Casson-Gordon invariants.
  \item $K$ has vanishing Ozsv\'ath-Szab\'o $\tau$-invariant and
    Rasmussen $s$-invariant.
  \end{enumerate}
\end{theorem}

We remark that in the proof of Theorem~\ref{theorem:main-example}~(1),
$g_*^t(K)$ is detected by considering a minimal genus problem for
which the results in~\cite{Hsiang-Szczarba:1971-1, Rokhlin:1971-1,
  Lee-Wilczynski:1997-1, Lee-Wilczynski:2000-1} give no interesting
lower bound but the $\rho$-invariants give an optimal bound.  We also
remark that results of Cochran-Orr-Teichner
\cite{Cochran-Orr-Teichner:1999-1} can be used to reveal partial
information that $g_*^t(K) > 0$, i.e., $K$ is not topologically slice.

As a consequence of Theorem~\ref{theorem:main-example}~(4), it follows
that the applications of the adjunction inequality to the smooth slice
genus as in \cite{Rudolph:1997-1, Lisca-Matic:1997-1,
  Lisca-Matic:1998-1, Akbulut-Matveyev:1997-1} give us no information
on~$K$, since $\tau$- and $s$-invariants are known to be sharper than
the Thurston-Bennequin lower bound, due to Plamenevskaya
\cite{Plamenevskaya:2004-1,Plamenevskaya:2004-2} and
Shumakovitch~\cite{Shumakovitch:2004-1}.  The author knows no other
method to apply gauge theory to estimate the slice genus of our~$K$.
Finally we remark that in the proof of
Theorem~\ref{theorem:main-example}~(4), we show a little more
generalized statement (Lemma~\ref{lemma:main-example}) that for any
finitely collection $\{\Phi_\alpha\}$ of integer-valued homomorphisms
of the smooth knot concordance group that give lower bounds of
$g_*^s$, our $K$ can be chosen in such a way that $\Phi_\alpha(K)=0$
for each~$\Phi_\alpha$, i.e., no such homomorphism extracts any
information on the slice genus of~$K$.  For more detailed discussion
on Theorem~\ref{theorem:main-example}, see
Section~\ref{section:slice-genus-of-knots}.

\section{Betti numbers and $L^2$-signatures}
\label{section:Betti-numbers-and-L2-signatures}

In this section we prove
Theorem~\ref{theorem:lower-bound-of-minimal-second-betti-number}~(1).
The essential part of the proof is to estimate the $L^2$-Betti number
of a 4-manifold in terms of the ordinary Betti number.  From this the
desired relationship between the ordinary Betti number and the
$L^2$-signature follows, because $L^2$-dimension theory enables us to
show that the $L^2$-signature is bounded by the (middle dimensional)
$L^2$-Betti number; this is an $L^2$-analogue of a well-known fact
that the ordinary signature is bounded by the Betti number.  In this
section all manifolds are topological manifolds.

\subsection*{Upper bounds of $L^2$-Betti numbers}

We start by defining the algebraic $L^2$-Betti number.  As a primary
reference on the $L^2$-theory we need, we refer to L\"uck's
book~\cite{Lueck:2002-1}.  Let $\Gamma$ be a discrete countable group.
While $L^2$-invariants are usually defined via the \emph{group von
  Neumann algebra $\N\Gamma$}, in this paper we will mainly use the
\emph{algebra $\U\Gamma$ of operators affiliated to~$\N\Gamma$}, which
is more useful for our purpose.  Both coefficients are known to give
the same $L^2$-Betti number and signature.

The $L^2$-dimension theory provides a dimension function
\[
\dim^{(2)}_\Gamma \colon \{ \text{finitely generated
  $\U\Gamma$-modules} \} \to [0,\infty).
\]
For a finite CW-complex $X$ endowed with $\pi_1(X) \to \Gamma$, the
twisted homology module
\[
H_i(X;\U\Gamma)=H_i(C_*(X;\Z\Gamma)\otimes_{\Z\Gamma} \U\Gamma)
\]
is defined by viewing $\U\Gamma$ as a $\Q\Gamma$-module via the
natural inclusions $\Q\Gamma \to \N\Gamma \to \U\Gamma$, and is known
to be finitely generated.  The \emph{$L^2$-Betti number}
$\beta_i^{(2)}(X)$ is defined to be $\beta_i^{(2)}(X) =
\dim^{(2)}_\Gamma H_i(X;\U\Gamma)$.  For a CW-pair $(X,A)$,
$\beta^{(2)}_i(X,A)$ is similarly defined.  (In this paper the choice
of $\pi_1(X)\to \Gamma$ will always be clearly understood and so we do
not include it in the notation.)  It is known that the analytic and
$L^2$-homological definitions are equivalent to the algebraic
definition described here \cite[Chapter 1, 6 and 8]{Lueck:2002-1}.

Following Cochran-Orr-Teichner~\cite{Cochran-Orr-Teichner:1999-1}, we
will focus on the case of a \emph{poly-torsion-free-abelian (PTFA)
  group}, which is defined to be a group admitting a finite length
normal series $\{G_i\}$ with torsion-free abelian quotients
$G_i/G_{i+1}$.  In this paper $\Gamma$ is always assumed to be PTFA.
Also, we assume that $\pi_1(X) \to \Gamma$ is nontrivial, since a
trivial homomorphism gives nothing beyond the (untwisted) rational
coefficient.

\begin{proposition}
  \label{proposition:upper-bound-of-L2-betti-number-of-4-manifold}
  Suppose $W$ is a connected compact 4-manifold (possibly with
  nonempty boundary) endowed with a nontrivial homomorphism $\pi_1(W)
  \to \Gamma$.  Then
  \begin{enumerate}
  \item $\beta_1^{(2)}(W) \le \beta_1(W)-1$,
  \item $\beta_2^{(2)}(W) \le \beta_2(W)$, and
  \item $\beta_3^{(2)}(W) \le 
    \begin{cases}
      \beta_3(W)-1 &\text{if $W$ is closed}, \\
      \beta_3(W) &\text{otherwise}.
    \end{cases}
    $
  \end{enumerate}
\end{proposition}

\begin{remark}
  \begin{enumerate}
  \item When $\partial W$ is nonempty, the proposition also gives an
    upper bound of $\beta^{(2)}_i(W,\partial W)$ in terms of the
    ordinary Betti number, by duality.
  \item In the special case that $H_1(\partial W;\Q) \cong H_1(W;\Q)$
    and $\partial W$ is nonempty, a similar result was proved (at
    least implicitly) in~\cite{Cochran-Orr-Teichner:1999-1}.  Our
    proof of
    Proposition~\ref{proposition:upper-bound-of-L2-betti-number-of-4-manifold}
    proceeds similarly to~\cite{Cochran-Orr-Teichner:1999-1}, but we
    need some technical modification to get rid of the
    $H_1$-isomorphism condition.
  \end{enumerate}
\end{remark}

Lemma~\ref{lemma:properties-of-PTFA-groups} below provides facts on a
PTFA group which are necessary to prove
Proposition~\ref{proposition:upper-bound-of-L2-betti-number-of-4-manifold}.
For a proof of Lemma~\ref{lemma:properties-of-PTFA-groups},
see~\cite{Cochran-Orr-Teichner:1999-1}.
 
\begin{lemma}
  \label{lemma:properties-of-PTFA-groups}
  \begin{enumerate}
  \item 
    $\Q\Gamma$ is an Ore domain so that there is a (skew-)field $\K$
    of right quotients of $\Q\Gamma$.  Every $\K$-module $M$ is free
    and has a well-defined dimension $\dim_\K M$.
    
  \item Suppose that $C_*$ is a finitely generated free chain complex
    over $\Q\Gamma$.  If $H_i(C_*\otimes_{\Q\Gamma} \Q)=0$ for $i \le
    n$, then $H_i(C_*\otimes_{\Q\Gamma} \K)=0$ for $i \le n$.
    
  \end{enumerate}
\end{lemma}

In particular, the existence of the skew-field $\K$ of quotients
enables us to understand the $L^2$-dimension as the ordinary dimension
over~$\K$, as follows: it is known that if $\Q\Gamma$ is an Ore
domain, then the natural map $\Q\Gamma \to \U\Gamma$ extends to an
embedding $\K\to \U\Gamma$~\cite{Lueck:2002-1}.  For a space $X$
equipped with $\pi_1(X) \to \Gamma$, let denote the Betti number with
$\K$-coefficients by $\beta_i(X;\K) = \dim_\K H_i(X;\K)$.  By
definition, $H_i(X;\U\Gamma)$ is the homology of the cellular chain
complex
\[
C_*(X;\Q\Gamma) \otimes_{\Q\Gamma} \U\Gamma =
(C_*(X;\Q\Gamma)\otimes_{\Q\Gamma} \K)\otimes_\K \U\Gamma.
\]
Since $H_*(X;\K)=H_*(C_*(X;\Q\Gamma)\otimes_{\Q\Gamma} \K)$, we have
the universal coefficient spectral sequence
\[
E_{p,q}^2 = \operatorname{Tor}^\K_p(H_q(X;\K),\U\Gamma) \Rightarrow
H_{p+q}(X;\U\Gamma).
\]
Since all higher $\operatorname{Tor}$ terms vanish over the
$\K$-coefficient, it follows that
\[
H_i(X;\U\Gamma)=H_i(X;\K) \otimes_\K \U\Gamma.
\]
Therefore $H_i(X;\U\Gamma)$ is always a free $\U\Gamma$-module whose
$\U\Gamma$-rank is equal to the $\K$-coefficient Betti number
$\beta_i(X;\K)$.  Since $\dim^{(2)}_\Gamma (\U\Gamma)^n=n$ (e.g.,
see~\cite{Lueck:2002-1}), we obtain

\begin{lemma}
  \label{lemma:L2-and-K-Betti-number}
  $\beta^{(2)}_i(X) = \beta_i(X;\K)$, and similarly for a pair
  $(X,A)$.
\end{lemma}

In order to prove
Proposition~\ref{proposition:upper-bound-of-L2-betti-number-of-4-manifold},
we first deal with the first Betti number.

\begin{lemma}
  \label{lemma:upper-bound-of-first-K-Betti-number}
  Suppose $(X,A)$ is a finite CW-pair with $X$ connected, and
  $\pi_1(X) \to \Gamma$ is a homomorphism.  Then
  \begin{enumerate}
  \item If $A$ is nonempty, $\beta_1(X,A;\K) \le \beta_1(X,A)$.
  \item If $A$ is empty, $\beta_1(X;\K) \le \beta_1(X)-1$.
  \end{enumerate}
\end{lemma}
We remark that the absolute case (2) was shown in \cite[Proposition
2.11]{Cochran-Orr-Teichner:1999-1}.

\begin{proof}
  Suppose that $A$ is nonempty.  Denote $\beta=\beta_1(X,A)$, and let
  $(Y,B)$ be the disjoint union of $\beta$ copies of $(I,\partial I)$
  where $I=[0,1]$.  Choose a map $f\colon (Y,B) \to (X,A)$ which
  induced an isomorphism $H_1(Y,B;\Q) \to H_1(X,A;\Q)$.

  By replacing $X$ with the mapping cylinder $M_f=(Y\times I) \cup
  X/(y,0)\sim f(x)$ of $f$, and replacing $A$ with $(B\times I) \cup A
  \subset M_f$, we may assume that $f$ is an injection $(Y,B)\subset
  (X,A)$ and $Y\cap A=B$.  From the homology long exact sequence with
  $\Q$-coefficients derived from
  \[
  0\to C_*(Y,B) \to C_*(X,A) \to C_*(X, Y\cup A) \to 0,
  \]
  it follows that $H_i(X, Y\cup A;\Q)=0$ for $i\le 1$.  By Lemma
  \ref{lemma:properties-of-PTFA-groups}~(2), $H_i(X,Y\cup A;\K)=0$ for
  $i\le 1$.  Thus, from the long exact sequence with
  $\K$-coefficients, it follows that $f$ induces a surjection
  $H_1(Y,B;\K) \to H_1(X,A;\K)$.  This shows that $\beta_1(X,A;\K) \le
  \beta_1(Y,B;\K)$.  On the other hand, since $C_i(Y,B;\K)=0$ for all
  $i$ but $C_1(Y,B;\K)=\K^\beta$, $\beta_1(Y,B;\K)= \beta$.  This
  completes the proof of~(1).

  Suppose $A$ is empty.  To apply the previous case, we choose a point
  $* \in X$ and consider the pair $(X,\{*\})$.  In the exact sequence
  \[
  0 \to H_1(X;\K) \to H_1(X,\{*\};\K) \to H_0(\{*\};\K) \to H_1(X;\K),
  \]
  $H_0(\{*\};\K)=\K$ obviously and $H_0(X;\K)=\K/(\pi_1(X)$-action) is
  trivial since $\K$ is a division ring and $\pi_1(X) \to \Gamma$ is
  nontrivial.  It follows that
  \[
  \beta_1(X;\K)+1 = \beta_1(X,\{*\};\K) \le \beta_1(X,\{*\}) =
  \beta_1(X).  \qedhere
  \]
\end{proof}

\begin{proof}[Proof of Proposition
  \ref{proposition:upper-bound-of-L2-betti-number-of-4-manifold}]
  Suppose $W$ is a compact connected 4-manifold equipped with
  $\pi_1(W) \to \Gamma$.  Since $W$ has the homotopy type of a finite
  CW-complex with cells of dimension $\le 4$, we may assume that the
  chain complex $C_*(W;-)$ is finitely generated and has dimension
  $\le 4$.

  By Lemma~\ref{lemma:L2-and-K-Betti-number}, we can think of
  $\beta_i(W;\K)$ instead $\beta^{(2)}_i(W)$.  So (1) follows directly
  from Lemma~\ref{lemma:upper-bound-of-first-K-Betti-number}.

  To prove (3), observe that the duality implies $\beta_3(W;\K)=
  \beta_1(W,\partial W;\K)$.  If $\partial W$ is empty,
  $\beta_1(W,\partial W;\K)\le \beta_1(W)-1 = \beta_3(W)-1$ by
  Lemma~\ref{lemma:upper-bound-of-first-K-Betti-number}.  If $\partial
  W$ is nonempty, $\beta_1(W,\partial W;\K)\le \beta_1(W,\partial W) =
  \beta_3(W)$ again by
  Lemma~\ref{lemma:upper-bound-of-first-K-Betti-number}.
  
  To prove (2), we use the fact that the Euler characteristics for the
  $\Q$- and $\K$-coefficients are the same, that is,
  \[
  \sum_{i=0}^4 (-1)^i \beta_i(W;\K) = \sum_{i=0}^4 (-1)^i \beta_i(W).
  \]
  Since $\pi_1(W) \to \Gamma$ is nontrivial, $\beta_0(W;\K)=0$.  When
  $W$ has nonempty boundary, $\beta_0(W,\partial W;\K) = 0$ since
  $\beta_0(W,\partial W;\K) \le \beta_0(W;\K)$.  From this it follows
  that $\beta_4(W;\K)=0$.  Plugging these values and the inequalities
  proved above into the Euler characteristic identity, we obtain
  $\beta_2(W;\K)\le \beta_2(W)$.
\end{proof}

\subsection*{Upper bounds of $L^2$-signatures}

We define the \emph{von Neumann $L^2$-signature} as follows: for a
$4k$-manifold $W$ endowed with $\pi_1(W) \to \Gamma$, the
$\U\Gamma$-coefficient intersection form
\[
\lambda\colon H_{2k}(W;\U\Gamma) \times H_{2k}(W;\U\Gamma) \to
\U\Gamma
\]
is a hermitian form.  In our case, $H_{2k}(W;\U\Gamma)$ is always a
free $\U\Gamma$-module since $\Gamma$ is assumed to be PTFA.  By
spectral theory, $H_{2k}(W;\U\Gamma)$ is decomposed as an orthogonal
sum of canonically defined subspaces $H_+$, $H_-$, and $H_0$ such that
$\lambda$ is positive definite, negative definite, and trivial, on
$H_+$, $H_-$, and $H_0$, respectively.  The $L^2$-signature of $W$ is
defined to be
\[
\sign^{(2)}(W)=\dim^{(2)}_\Gamma(H_+)-\dim^{(2)}_\Gamma(H_-) \in \R.
\]
For more details and the relationship with other ways to define the
$L^2$-signature, refer to~\cite{Cochran-Orr-Teichner:1999-1}
and~\cite{Lueck-Schick:2003-1}.
 
\begin{lemma}
  \label{lemma:L2-signature-and-L2-Betti-number}
  $|\sign^{(2)}(W)| \le \beta^{(2)}_{2k}(W)$.
\end{lemma}
\begin{proof}
  Since $H_+, H_- \subset H_{2k}(W;\U\Gamma)$ and $H_+\cap H_-=\{0\}$,
  $L^2$-dimension theory enables us to show
  \[
  \dim^{(2)}_\Gamma(H_+)+\dim^{(2)}_\Gamma(H_-) \le \dim^{(2)}_\Gamma
  H_{2k}(W;\U\Gamma)
  \]
  using an $L^2$-analogue of a standard argument of elementary linear
  algebra.  (e.g., refer to Chapter 8 of \cite{Lueck:2002-1}, where it
  is shown that $\dim^{(2)}_\Gamma$ satisfies a set of axioms which
  includes all the properties we need.)  From this the conclusion
  follows.
\end{proof}

Combining Lemma~\ref{lemma:L2-signature-and-L2-Betti-number} with
Lemma~\ref{proposition:upper-bound-of-L2-betti-number-of-4-manifold},
we obtain:

\begin{lemma}
  \label{lemma:L2-signature-and-Betti-number}
  If $W$ is a compact connected 4-manifold endowed with a nontrivial
  homomorphism $\pi_1(W) \to \Gamma$, then
  \[
  |\sign^{(2)}(W)| \le \beta_2(W).
  \]
\end{lemma}

Now we are ready to show the first part of
Theorem~\ref{theorem:lower-bound-of-minimal-second-betti-number}
stated in the introduction.  We adopt the following topological
definition of the $\rho$-invariant, as in
Chang-Weinberger~\cite{Chang-Weinberger:2003-1}.  (See also
Cochran-Orr-Teichner~\cite{Cochran-Orr-Teichner:1999-1}.)  Let $M$ be
a 3-manifold endowed with $\pi_1(M) \to \Gamma$.  It is known that
there is a bigger group $G$ containing $\Gamma$ and a 4-manifold $W$
such that $\partial W$ consists of $r$ components $M_1, \ldots, M_r$
($r>0$), $M_i \cong M$, and $\pi_1(M_i) \xrightarrow{\phi} \Gamma \to
G$ factors through $\pi_1(W)$ for each~$i$.  (For a proof, see the
appendix of~\cite{Chang-Weinberger:2003-1}; they consider the case
that $\pi_1(M)=\Gamma$ but the same argument works in our case as
well.)  Then $\rho(M,\phi)$ is defined to be the following signature
defect:
\[
\rho(M,\phi) = \frac{1}{r} \big( \sign^{(2)}(W) - \sign(W) \big) \in
\R
\]
where $\sign^{(2)}(W)$ and $\sign(W)$ denote the $L^2$-signature
associated to $\pi_1(W)\to G$ and the ordinary signature,
respectively.  The real number $\rho(M,\phi)$ is determined by $M$ and
$\phi$, and independent of the choices we made.  From the results
in~\cite{Lueck-Schick:2003-1} it follows that $\rho(M,\phi)$ defined
above coincides with the $\rho$-invariant of
Cheeger-Gromov~\cite{Cheeger-Gromov:1985-1}.

\begin{proof}[Proof of
  Theorem~\ref{theorem:lower-bound-of-minimal-second-betti-number}~(1)]
  Suppose $W$ is a compact connected 4-manifold with boundary $M$, and
  $\pi_1(W) \to \Gamma$ is given.  Let denote the composition
  $\pi_1(M) \to \pi_1(W) \to \Gamma$ by~$\phi$.  Our goal is to show
  that $|\rho(M,\phi)| \le 2\beta_2(W)$.

  Since $\phi$ factors through $\pi_1(W)$, we can compute
  $\rho(M,\phi)$ using~$W$; by the definition above,
  \[
  \rho(M,\phi) = \sign^{(2)}(W)-\sign(W).
  \]
  Obviously $|\sign(W)| \le \beta_2(W)$.  By
  Lemma~\ref{lemma:L2-signature-and-Betti-number}, $|\sign^{(2)}(W)|
  \le \beta_2(W)$.  From this the desired conclusion follows.
\end{proof}

\section{Extending coefficient systems to bounding 4-manifolds}
\label{section:extending-coefficient-systems}

Suppose $W$ is a topological 4-manifold with boundary $M$ and
$\pi_1(W) \to \Gamma$ is given.  ($M$~is endowed with the induced map
$\pi_1(M) \to \Gamma$.)  In this section we deal with the problem of
extending a bigger coefficient system on $M$ to $W$ to prove
Theorem~\ref{theorem:lower-bound-of-minimal-second-betti-number}~(2).
To state a more detailed form of
Theorem~\ref{theorem:lower-bound-of-minimal-second-betti-number}~(2),
we need the following facts from~\cite{Cochran-Orr-Teichner:1999-1}:
suppose $\sR$ is a (possibly non-commutative) subring of $\K$ which is
a PID containing~$\Z\Gamma$.  In this section we assume that
$H_1(M;\sR)$ is $\sR$-torsion.

\medskip

(1) \emph{Blanchfield form on $H_1(M;\sR)$.}  The Bockstein map $B\colon
H_2(M;\K/\sR) \to H_1(M;\sR)$ and the Kronecker evaluation $\kappa\colon
H^1(M;\K/\sR) \to \Hom(H_1(M;\sR), \K/\sR)$ are isomorphisms.  The
Blanchfield form, which is defined to be the isomorphism
\begin{multline*}
  B\ell\colon H_1(M;\sR) \xrightarrow{B^{-1}} H_2(M;\K/\sR)
  \xrightarrow{\text{duality}} H^1(M;\K/\sR) \\ \xrightarrow{\kappa}
  \Hom(H_1(M;\sR),\K/\sR),
\end{multline*}
is a symmetric linking form on $H_1(M;\sR)$ \cite[p.\
451]{Cochran-Orr-Teichner:1999-1}.

\medskip

(2) \emph{Coefficient systems induced by characters.} A homomorphism
$h\colon H_1(M;\sR) \to \K/\sR$ gives rise to a group homomorphism
$\psi\colon \pi_1(M) \to \K/\sR\rtimes \Gamma$ in a natural way.
Indeed, $\psi$ is a lift of $\pi_1(M) \to \Gamma$, i.e.,
\[
\begin{diagram}
  \node[2]{\K/\sR \rtimes \Gamma} \arrow{s,r}{\text{projection}}
  \\
  \node{\pi_1(M)} \arrow{e} \arrow{ne,t,..}{\psi} \node{\Gamma}
\end{diagram}
\]
commutes, and the restriction of $\psi$ on $N=\Ker\{\pi_1(M) \to
\Gamma\}$ agrees with
\[
N \to N/[N,N]=H_1(M;\Z\Gamma) \to H_1(M;\sR) \xrightarrow{h} \K/\sR
\subset \K/\sR \rtimes \Gamma.
\]
Furthermore, $\psi$ factors through $\pi_1(W)$ if $h$ factors through
$H_1(W;\sR)$ \cite[p.\ 455]{Cochran-Orr-Teichner:1999-1}.

\medskip

Note that $\K/\sR$ is a torsion-free abelian group, and therefore
$\K/\sR \rtimes \Gamma$ is PTFA when $\Gamma$ is PTFA.  We also recall
that, as in case of a commutative PID, any finitely generated
$\sR$-module $M$ is isomorphic to $F\oplus tM$ where $F$ is a free
module of rank $\dim_\K (M\otimes_\sR \K)$ and $tM$ is the
$\sR$-torsion submodule of~$M$.  (e.g., refer to~\cite{Cohn:1971-1}.)
$tM$ is isomorphic to a direct sum of cyclic modules of nonzero order.

Now we can state the result we will prove in this section.  Denote by
$\partial$ the boundary map $H_2(W,M;\sR) \to H_1(M;\sR)$.

\begin{theorem}
  \label{theorem:character-extension}
  Suppose that $H_2(W,M;\sR)=F\oplus tH_2(W,M;\sR)$ and $\partial(F)$ is a
  proper submodule of $H_1(M;\sR)$ for some free summand~$F$.  Then
  there is a nontrivial submodule $P$ in $H_1(M;\sR)$ such that for any
  $x\in P$, the homomorphism
  \[
  B\ell (x) \colon H_1(M;\sR) \to \K/\sR
  \]
  factors through $H_1(W;\sR)$.
\end{theorem}

In particular, if $H_1(M;\sR)$ is never generated by $\beta_2(W)$
elements, then since
\[
\beta_2(W,M;\K) = \beta_2(W;\K) = \beta^{(2)}_2(W) \le \beta_2(W)
\]
by duality and
Lemma~\ref{proposition:upper-bound-of-L2-betti-number-of-4-manifold},
the hypothesis of Theorem~\ref{theorem:character-extension} is
satisfied.  It it is the case, then for $x \in P$, $B\ell(x)$ gives
rise to a homomorphism $\pi_1(M) \to \K/\sR \rtimes \Gamma$ which
factors through~$\pi_1(W)$.  This proves
Theorem~\ref{theorem:lower-bound-of-minimal-second-betti-number}~(2).

The remaining part of this section is devoted to the proof of
Theorem~\ref{theorem:character-extension}.  As the first step, we will
show that for the boundary of a relative 2-cycle of $(W,M)$, the
Blanchfield form of $M$ can be computed via the intersection form
of~$W$.  Indeed it is a consequence of the following algebraic
observation:

\begin{lemma}
  \label{lemma:bockstein-and-pair}
  Suppose $\sR$ is a (possibly non-commutative) ring with
  (skew-)quotient field $\K$, and
  \[
  0 \to C'_* \xrightarrow{i} C_* \xrightarrow{p} C''_* \to 0
  \]
  is an exact sequence of chain complexes over $\sR$ such that
  $H_n(C'\otimes \K)=0=H_{n-1}(C'\otimes \K)$.  Then
\long\def\ignoreme{
the diagram
  \[
  \begin{diagram}\dgARROWLENGTH=1.5em \dgHORIZPAD=0mm \dgVERTPAD=0mm
    \node[2]{\b{H_n(C_*)}} \arrow{e,t}{p_*} \arrow{s}
    \node{\b{H_n(C''_*)}} \arrow{s} \arrow{se,-}
    \\
    \node[2]{\b{H_n(C_*\otimes\K)}} \arrow{e,tb}{p_*}{\cong} \arrow{s}
    \node{\b{H_n(C''_*\otimes\K)}} \arrow{s}
    \node{} \arrow{s,r,-}{\partial}
    \\
    \node{\b{H_n(C'_*\otimes\K/\sR)}} \arrow{s,lr}{\cong}{B=\text{Bockstein}}
      \arrow{e,t}{i_*}
    \node{\b{H_n(C_*\otimes \K/\sR)}} \arrow{e,t}{p_*}
    \node{\b{H_n(C''_*\otimes \K/\sR)}}
    \node{} \arrow{sw,-}
    \\
    \node{\b{H_{n-1}(C'_*)}} \node[2]{} \arrow[2]{w}
  \end{diagram}    
  \]
  commutes, i.e.,
}
  \[
  \alpha\colon H_n(C''_*) \to H_n(C''_* \otimes \K) \xrightarrow[\cong]{p_*^{-1}}
  H_n(C_* \otimes \K) \to H_n(C_* \otimes \K/\sR)
  \]
  coincides with
  \[
  \beta\colon H_n(C''_*) \xrightarrow{\partial} H_{n-1}(C'_*)
  \xrightarrow[\cong]{B^{-1}} H_n(C'_* \otimes \K/\sR) \xrightarrow{i_*}
  H_n(C\otimes \K/\sR).
  \]
\end{lemma}

\begin{proof}
  First note that the Bockstein $B$ and the induced map $p_*$ are
  isomorphisms since $H_n(C'_*\otimes \K)=0=H_{n-1}(C'_*\otimes \K)$.

  We will regard $C'_*$ as a submodule of $C_*$ and denote the
  homology class of a cycle $x$ by~$[x]$.  Suppose $z$ is a cycle
  in~$C''_n$, and $x\in C_n$ is a pre-image of $z$, i.e., $p(x)=z$.
  $H_j(C'_*)\otimes \K=H_j(C'_* \otimes \K)=0$ for $j=n$, $n-1$ since
  $\K$ is a flat $\sR$-module, and therefore $p$ induces an isomorphism
  $H_n(C_*)\otimes \K \cong H_n(C''_*)\otimes \K$.  It follows that
  there is a cycle $y$ in $C_n$ such that $p_*[y]=[z]\cdot r$ in
  $H_n(C''_*)$ for some nonzero $r\in \sR$, that is, there is $u\in
  C_{n+1}$ such that $\partial u = x\cdot r-y+w$ where $w\in C'_n
  \subset C_n$.  Since $\partial y=0$, $\partial w = \partial x \cdot
  r$.  Therefore $w \otimes \frac{1}{r}$ is a cycle in $C'_n \otimes
  \K/\sR$.

  Since $p_*[y\otimes \frac 1r]=[z\otimes 1]$, it can be seen that
  $\alpha[z]=[y \otimes \frac{1}{r}]$.  On the other hand, by the
  definition of the Bockstein homomorphism, $B[w\otimes \frac
  1r]=[\partial x]$, and therefore, $\beta[z]=[w\otimes \frac 1r]$.

  In $C_n \otimes \K/\sR$, we have
  \[
  \partial \Big(u\otimes \frac 1r\Big) = x-y\otimes \frac 1r + w \otimes \frac
  1r = -y\otimes \frac 1r + w\otimes \frac 1r.
  \]
  From this it follows that $[y\otimes \frac 1r] = [w\otimes \frac
  1r]$ in $H_n(C_*\otimes \K/\sR)$.
\end{proof}

Recall that $\partial$ denotes the boundary map $H_2(W,M;\sR) \to
H_1(M;\sR)$.

\begin{lemma}
  \label{lemma:Blanchfield-form-for-boundary}
  Let $\Phi$ be the composition
  \begin{multline*}
    \Phi\colon H_2(W,M;\sR) \to H_2(W,M;\K)
    \cong H_2(W;\K) \to H_2(W;\K/\sR)\\
    \cong H^2(W,M;\K/\sR) \xrightarrow{\kappa}
    \Hom(H_2(W,M;\sR),\K/\sR),
  \end{multline*}
  where $\kappa$ is the Kronecker evaluation map.  Then
  $B\ell(\partial x)(\partial y) = \Phi(x)(y)$ for any $x, y \in
  H_2(W,M;\sR)$.
\end{lemma}

\begin{proof}
  From Lemma~\ref{lemma:bockstein-and-pair} and the naturality of
  duality and the Kronecker evaluation, we obtain the following
  commutative diagram:
  \[
  \begin{diagram}\dgARROWLENGTH=1.7em \dgHORIZPAD=0mm \dgVERTPAD=0mm
    \def\b#1{\fbox{$\displaystyle #1$}} \fboxsep=.5em \fboxrule=0pt
    \node[2]{\b{H_2(W,M;\sR)}} \arrow{e,t}{\partial} \arrow{s}
    \node{\b{H_1(M;\sR)}} \arrow{se,-}
    \\
    \node{\b{H_2(W;\K)}} \arrow{e,b}{\cong} \arrow{s}
    \node{\b{H_2(W,M;\K)}} \node[2]{}\arrow[2]{s,r,-}{B\ell}
    \\
    \node{\b{H_2(W;\K/\sR)}} \arrow{s,lr}{\mathrm{duality}}{\cong}
    \node[2]{\b{H_1(M;\K/\sR)}}
    \arrow[2]{n,lr}{\cong}{B} \arrow[2]{w}
    \arrow{s,lr}{\cong}{\mathrm{duality}}
    \\
    \node{\b{H^2(W,M;\K/\sR)}} \arrow{s,l}{\kappa}
    \node[2]{\b{H^1(M;\K/\sR)}} \arrow[2]{w} \arrow{s,lr}{\cong}{\kappa}
    \node{} \arrow{sw}
    \\
    \node{\b{\Hom(H_2(W,M;\sR),\K/\sR)}}
    \node[2]{\b{\Hom(H_1(M;\sR),\K/\sR)}} \arrow[2]{w,t}{\partial^\#}
  \end{diagram}    
  \]
  Here the map $\partial^\#$ is given by
  $\partial^\#(\psi)(y)=\psi(\partial y))$ for $\psi\colon H_1(M;\sR)
  \to \K/\sR$ and $y \in H_2(W,M;\sR)$.  From this the conclusion
  follows.
\end{proof}

\begin{proof}[Proof of Theorem~\ref{theorem:character-extension}]
  For a submodule $P$ in $H_1(M;\sR)$, we denote
  \[
  P^{\perp} = \{ y\in H_1(M;\sR) \mid B\ell(x)(y)=0 \text{ for all }x\in P \}.
  \]
  Consider the exact sequence
  \[
  \cdots \to H_2(W;\sR) \to H_2(W,M;\sR) \xrightarrow{\partial} H_1(M;\sR)
  \to H_1(W;\sR) \to \cdots.
  \]
  We will show that there is a nontrivial submodule $P$ in $H_1(M;\sR)$
  such that the image $\partial(H_2(W,M;\sR))$ is contained
  in~$P^\perp$.  Indeed from this claim it follows that, for any $x\in
  P$, $B\ell(x)\colon H_1(M;\sR) \to \K/\sR$ gives rise to a homomorphism
  $\operatorname{Coker} \partial \to \K/\sR$, which automatically
  extends to $H_1(W;\sR)$ since $\K/\sR$ is an injective $\sR$-module.  This
  completes the proof.

  Recall that we wrote $H_2(W,M;\sR)=F \oplus tH_2(W,M;\sR)$ where $F$ is
  free and $tH_2(W,M;\sR)$ is the torsion submodule.  To prove the
  claim, we consider the following two cases:

  Case 1: Suppose $\partial(tH_2(W,M;\sR))$ is nontrivial.  Consider
  the composition $\Phi$ described in
  Lemma~\ref{lemma:Blanchfield-form-for-boundary}.  For any $x\in
  tH_2(W,M;\sR)$ we have $\Phi(x)=0$, since $\Phi$ factors through
  $H_2(W,M;\K)$ which is torsion free.  Therefore $B\ell(\partial
  x)(\partial y)=\Phi(x)(y)=0$ for any $y\in H_2(W,M;\sR)$.  This
  shows that $P=\partial(tH_2(W,M;\sR))$ is a nontrivial submodule
  satisfying the desired property.

  Case 2: Suppose $\partial(tH_2(W,M;\sR))$ is trivial.  Then the image
  of $\partial$ is equal to $\partial(F)$, which is a proper submodule
  of $H_1(M;\sR)$ by the hypothesis.  Appealing to the lemma below,
  which should be regarded as folklore, it follows that $P=\partial
  (F)^\perp$ is nontrivial.  It is obvious that this $P$ has the
  desired property.
\end{proof}

\begin{lemma}
  Suppose $A$ is a finitely generated torsion $\sR$-module endowed
  with a symmetric linking form given by an isomorphism $\Psi \colon A
  \to \Hom(A,\K/\sR)$.  Then for any proper submodule $B$ in $A$,
  $B^\perp$ is nontrivial.
\end{lemma}

\begin{proof}
  From the exact sequence
  \[
  0 \to \Hom(A/B,\K/\sR) \xrightarrow{p^\#} \Hom(A,\K/\sR) \xrightarrow{i^\#}
  \Hom(B,\K/\sR)
  \]
  it follows that $B^\perp = \Psi^{-1}(\Ker i^\#) = \Psi^{-1}(\Im
  p^\#)$.  So it suffices to show that $\Hom(A/B,\K/\sR)$ is
  nontrivial.  Note that every cyclic module $\sR/p\sR$ with $p\ne 0$
  is (isomorphic to) a submodule of~$\K/\sR$ .  Since $A/B$ is a
  nontrivial torsion module, it has a summand of the form $\sR/p\sR$
  with $p\ne 0$, by the structure theorem of finitely generated
  $\sR$-modules.  It follows that $\Hom(A/B,\K/\sR)$ is nontrivial.
\end{proof}

\section{Construction of a bordism from a locally flat surface}
\label{section:construction-of-bordisms-from-locally-flat-surfaces}

In this section we will prove
Proposition~\ref{proposition:Z-bordism-from-embedded-surface}.
Suppose $W$ is a topological 4-manifold with boundary $M$ such that
$H_1(W)=0$, and $\sigma$ is a 2-dimensional homology class contained
in $\Im\{H_2(M) \to H_2(W)\}$.  First we describe a homomorphism
$\phi_\sigma \colon \pi_1(M) \to \Z$ which is determined by~$\sigma$.
Consider the exact sequence
\[
H_2(W) \to H_2(W,M) \xrightarrow{\partial} H_1(M) \to H_1(W)=0.
\]
The intersection with~$\sigma$ gives a homomorphism $\sigma\cdot
\colon H_2(W,M) \to \Z$, which induces a homomorphism $h_\sigma \colon
H_1(M) \to \Z$ since $\sigma\cdot$ vanishes on the image of~$H_2(W)$.
Define $\phi_\sigma$ to be the composition
\[
\phi_\sigma\colon \pi_1(M) \to H_1(M) \xrightarrow{h_\sigma} \Z.
\]

Recall that
Proposition~\ref{proposition:Z-bordism-from-embedded-surface} claims
that if there is a locally flat surface $\Sigma$ of genus $g$ in $W$
which represents the class $\sigma \in H_2(W)$ and the map
$\phi_\sigma$ is nontrivial, then there is a topological 4-manifold
$V$ bounded by $M$ such that $\beta_2(V) = \beta_2(W)+2g-1$ and
$\phi_\sigma$ factors through~$H_1(V)$.  Roughly speaking, we will
construct $V$ by performing ``surgery along $\Sigma$'' on~$W$.

\begin{proof}[Proof of
  Proposition~\ref{proposition:Z-bordism-from-embedded-surface}]

  By Alexander duality, $H_2(W,W-\Sigma)$ can be identified with
  $H^2(\Sigma) = \Z$.  From the exact sequence
  \[
  H_2(W) \xrightarrow{\sigma\cdot} H_2(W,W-\Sigma) \to H_1(W-\Sigma)
  \to H_1(W)=0
  \]
  it follows that $H_1(W-\Sigma)\cong H_2(W,W-\Sigma)=\Z$ since the
  leftmost map $\sigma\cdot$ is given by the intersection of a 2-cycle
  with $\sigma$, which is always zero.

  Note that $\Sigma$ has trivial normal bundle in $W$ since $\Sigma$
  is connected and the self-intersection $\sigma\cdot \sigma$
  vanishes.  There is a bijection between the set of (fiber homotopy
  classes of) framings on $\Sigma$ and $[\Sigma,S^1]=H^1(\Sigma,\Z)$
  which can be identified with $\Z^{2g}$ by choosing a basis $\{x_i\}$
  of $H_1(\Sigma)$.  Pushoff along a framing induces a homomorphism
  $H_1(\Sigma) \to H_1(W-\Sigma)$ in such a way that if the framing
  corresponding to $0 \in \Z^{2g}$ induces $h\colon H_1(\Sigma) \to
  H_1(W-\Sigma)$, then the framing corresponding to $(a_i)\in \Z^{2g}$
  gives rise to a homomorphism sending $x_i$ to $h(x_i)+a_i[\mu]$
  where $\mu$ is a meridional curve of~$\Sigma$.  Since
  $H_1(W-M)\cong\Z$ is generated by $[\mu]$, it follows that there is
  a framing inducing a trivial homomorphism $H_1(\Sigma) \to
  H_1(W-\Sigma)$.  We identify a tubular neighborhood of $\Sigma$ in
  $W$ with $\Sigma\times D^2$ under this framing, and denote
  $N=W-\operatorname{int}(\Sigma\times D^2)$.

  Choose a 3-manifold $R$ with boundary $\Sigma$ such that
  $H_1(\Sigma) \to H_1(R)$ is surjective (e.g., a handlebody with the
  same genus as $\Sigma$ may be used as $R$).  Let
  \[
  V=\big(N\cup (R\times S^1)\big)/\sim
  \]
  where $\Sigma\times S^1 \subset \partial N$ and $\partial R \times
  S^1$ are identified.  From the Mayer-Vietoris sequence
  \[
  \cdots \to H_1(\Sigma\times S^1) \to H_1(N) \oplus H_1(R\times S^1)
  \to H_1(V) \to 0
  \]
  for $V=N\cup (R\times S^1)$, it follows that $H_1(V) \cong H_1(N) =
  \Z$ since $H_1(\Sigma) \to H_1(R)$ is surjective and $i_*\colon
  H_1(\Sigma) \to H_1(N)$ is trivial by our choice of the framing
  on~$\Sigma$.  From the definition it is easily seen that $h_\sigma$
  is equal to the map $H_1(M) \to H_1(V)=\Z$ induced by the inclusion.
  Therefore $\phi_\sigma$ factors through~$\pi_1(V)$ as desired.

  The Betti number assertion follows from a straightforward
  computation.  For the convenience of the reader, we give details
  below.  From the above Mayer-Vietoris sequence it follows that
  \[
  \chi(\Sigma\times S^1)+\chi(V)=\chi(N)+\chi(R\times S^1)
  \]
  where $\chi$ denotes the Euler characteristic.
  $\chi(N)+\chi(\Sigma)=\chi(W)$ by the long exact sequence for the
  pair $(W,N)$ and Alexander duality.  Since $\chi(X\times S^1)=0$ for
  any $X$, it follows that
  \[
  \chi(V)=\chi(W)-\chi(\Sigma)=\chi(W)+2g-2.
  \]
  From the hypothesis that $H_1(W)=0$, it follows that $\beta_1(W)=0$
  and $\beta_3(W)=\beta_1(W,M)=\beta_0(M)-1$.  $\beta_1(V)=1$ as shown
  above.  Since $\phi_\sigma$ is nontrivial, so is $H_1(M) \to
  H_1(V)=\Z$ and thus has torsion cokernel.  It follows that
  $\beta_3(V)=\beta_1(V,M)=\beta_0(M)-1$.  Combining these
  observations on the Betti numbers with the Euler characteristic
  identity, the desired inequality follows.
\end{proof}

\section{Slice genus}
\label{section:slice-genus-of-knots}

In this section we apply the results proved in the previous sections
to investigate the slice genus of a knot $K$ in~$S^3$.  Indeed our
results give lower bounds of the genus of a spanning surface in a
homology 4-ball; for a knot $K$ in a homology 3-sphere $Y$ which
bounds some (topological) homology 4-ball, let $g_*^h(K)$ be the
minimal genus of a locally flat surface $F$ in a homology 4-ball $X$
such that $\partial (X,F)=(Y,K)$.  Obviously $g_*^h(K) \le g_*^t(K)
\le g_*^s(K)$ for a knot $K$ in~$S^3$.

For $(X,F)$ as above, consider the 4-manifold $W$ obtained by
attaching a 2-handle to $X$ along the preferred framing of~$K$.  The
boundary of $W$ is the result of surgery on $Y$ along the preferred
framing of~$K$, which we will call the \emph{zero-surgery manifold of
  $K$} and denote by~$M_K$.  Note that $H_1(M_K)=\Z$ is generated by a
meridian of~$K$.  Let $\sigma$ be a generator of $H_2(W)=\Z$.  It can
be easily seen that the abelianization map $\phi\colon\pi_1(M_K) \to
H_1(M_K)=\Z$ is exactly the homomorphism $\phi_\sigma$ defined in
Section~\ref{section:construction-of-bordisms-from-locally-flat-surfaces}.
Also, note that $\sigma$ is represented by a surface in $M_K$, namely a
capped-off Seifert surface of~$K$.

Attaching to $F$ the core of the 2-handle of $W$, we obtain a surface
$\Sigma$ with the same genus as $F$ which represents the homology
class $\sigma\in H_2(W)$.  Therefore, by
Proposition~\ref{proposition:Z-bordism-from-embedded-surface}, one
obtains a null-bordism of $M_K$ over $\Z$ with bounded~$\beta_2$; we
state it as a proposition.

\begin{proposition}
  \label{proposition:bordism-from-slice-surface}
  There is a topological 4-manifold $V$ with boundary $M_K$ such that
  $\phi\colon \pi_1(M_K) \to \Z$ factors through~$\pi_1(V)$ and
  $\beta_2(V) \le 2g_*^h(K)$.  In particular, if $K$ is a knot in
  $S^3$, then $\beta_2(V) \le 2g^t_*(K)$.
\end{proposition}

This enables us to use
Theorem~\ref{theorem:lower-bound-of-minimal-second-betti-number},
possibly repeatedly, to obtain lower bounds of~$g^h_*(K)$.  We remark
that while a lower bound is obtained from $\rho(M_K,\phi)$ by applying
Theorem~\ref{theorem:lower-bound-of-minimal-second-betti-number}~(1)
directly, it gives us no interesting result since it is known that
$\rho(M_K,\phi)$ is determined by the signature function of $K$
\cite{Cochran-Orr-Teichner:1999-1}.  However, it turns out that the
$\rho$-invariants associated to bigger coefficient systems obtained by
Theorem~\ref{theorem:lower-bound-of-minimal-second-betti-number}~(2)
actually reveal new information on the slice genus which cannot be
obtained via previously known invariants, as mentioned in
Theorem~\ref{theorem:main-example}.  The remaining part of this
section is devoted to the construction of examples illustrating this.

\subsection*{Construction of examples}

Our examples will be constructed using a well known method that
produces a new knot from a given knot by ``tying'' another knot along
a circle in the complement.  For a knot $J$, we denote its exterior by
$E_J=S^3-$ (open tubular neighborhood of~$J$).  Suppose $K_0$ is a
knot and $\eta$ is a circle in $S^3-K_0$ which is unknotted in~$S^3$.
Choose a (closed) tubular neighborhood $U$ of~$\eta$ in $S^3-K_0$.
Removing the interior of $U$ from $S^3-K_0$ and attaching the exterior
$E_{J}$ of a knot $J$ along the boundary of $U$ in such a way that a
meridional curve of $\eta$ is identified with a curve null-homologous
in $E_{J}$, one obtains the complement of a new knot in $S^3$, which
we will denote by $K_0(\eta, J)$.  In some literature this
construction is called the ``satellite construction'' or ``genetic
infection''.

We start by choosing a knot $K_s$ in $S^3$ whose Alexander polynomial
$\Delta_{K_s}(t)$ is a cyclotomic polynomial $\Phi_n(t)$ with $n$
divisible by at least three distinct primes.  Indeed, by a well-known
characterization due to Levine, there is such a knot if and only if
$\Phi_n(t^{-1})=\pm t^s \Phi_n(t)$ for some $s$ and $\Phi_n(1)=\pm 1$.
Since the complex conjugate of a root of unity is its reciprocal,
$\Phi_n(t)$ satisfies the former condition.  For the latter condition,
one may appeal to the following lemma:

\begin{lemma}
  For $n\ge 2$, $\Phi_n(1)=1$ if and only if $n$ is not a prime power.
\end{lemma}

\begin{proof}
  If $n=p^a$ is a prime power, then it is easily seen that $\Phi_n(t)$ is
  given by
  \[
  \Phi_n(t) = t^{p^{a-1}(p-1)} + \cdots t^{p^{a-1}}+1
  \]
  and therefore $\Phi_n(1)=p$.

  Conversely, suppose $n=p_1^{a_1} \cdots p_r^{a_r}$ with $p_i$ prime
  and $r>1$.  We recall that
  \[
  t^n-1=\prod_{d \mid n} \Phi_d(t).
  \]
  By eliminating the factor of $t-1$ and rearranging terms, we obtain
  \[
  t^{n-1} + \cdots + t + 1 = \bigg( \prod_{i=1}^r \prod_{j=1}^{a_i}
  \Phi_{p_i^j}(t) \bigg) \cdot \Phi_n(t) \cdot h(t)
  \]
  Plugging $t=1$, it follows that $\Phi_n(1)h(1)=1$ and so
  $\Phi_n(1)=1$.
\end{proof}

Denote the (rational) Alexander module $H_1(M_{J};\Q[t,t^{-1}])$ of a
knot $J$ by~$A_{J}$, and the mirror image of $J$ by~$-J$.  (Here we
adopt the standard convention of the orientation of $-J$ so that $J\#
(-J)$ is always a ribbon knot.)

Returning to our construction, for an unknotted circle in $\eta$
disjoint to $K_s$ and two knots $J$ and $J'$ which will be chosen
later, consider the connected sum
\[
K = \mathop{\#}^g  \big(K_s(\eta, J) \# -(K_s(\eta,J')) \big)
\]
of $g$ identical knots.

We choose $\eta$ in such a way that the following properties are
satisfied:

\begin{enumerate}
\item[(P1)] The linking number of $\eta$ and $K_s$ vanishes, so that
  $\eta$ represents a homology class $[\eta]\in A_{K_s}$.
  Furthermore, $[\eta]$ is a generator of~$A_{K_s}$.
\item[(P2)] For any $J$ and $J'$, $K$ satisfies $g_*^s(K) \le g$.
\item[(P3)] For any $J$ and $J'$, $K$ is algebraically slice, i.e.,
  $K$ has a Seifert matrix of a slice knot.
\item[(P4)] For any $J$ and $J'$, $K$ has vanishing Casson-Gordon
  invariants.
\end{enumerate}

\begin{figure}[ht]
  \begin{center}
    \includegraphics{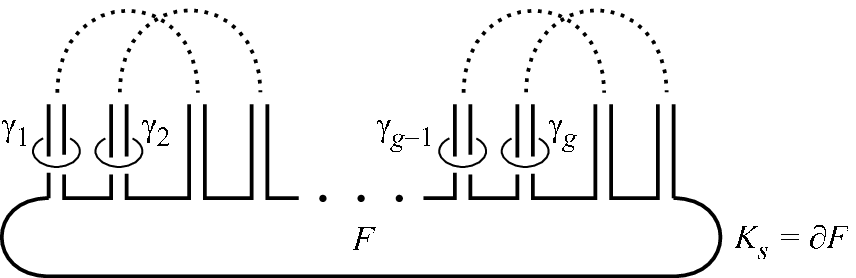}
  \end{center}
  \caption{}
  \label{figure:seed-knot}
\end{figure}

For this purpose, we first choose a Seifert surface $F$ of~$K_s$.
$F$~consists of one 0-handle and $2r$ 1-handles, where $r$ is the
genus of~$F$.  Choose unknotted circles $\gamma_1, \ldots,
\gamma_{2r}$ in $S^3-F$ which are Alexander dual to the 1-handles
of~$F$, as illustrated in Figure~\ref{figure:seed-knot}.

Since each $\gamma_i$ is disjoint to~$F$, it represents a homology
class $[\gamma_i] \in A_{K_s}$.  Also, it can be seen that the
$[\gamma_i]$ generate~$A_{K_s}$, by a standard Mayer-Vietoris
argument.  Therefore one of the $[\gamma_i]$, say $[\gamma_1]$, is
nontrivial in $A_{K_s}$.  Let $\eta$ be $\gamma_1$.

\begin{lemma}
  $\eta$ satisfies the properties (P1)--(P4) required above.
\end{lemma}

\begin{proof}
  Obviously $\eta$ has linking number zero with~$K$.  Since
  $\Delta_{K_s}(t)$ is irreducible, $A_{K_s}=\Q[t,t^{-1}]/\langle
  \Delta_{K_s}(t) \rangle$, and $[\eta]\ne 0$ is automatically a
  generator of~$A_{K_s}$.  This shows~(P1).

  Let $L = K_s(\eta, J) \# -(K_s(\eta,J'))$.  We claim that
  $g_*^s(L)\le 1$, from which (2) easily follows.  To prove the claim,
  observe that $L$ is obtained from the ribbon knot $K_s \# (-K_s)$,
  by ``tying'' $J$ and~$J'$.  Note that the boundary connected sum of
  $F$ and $-F$ is a Seifert surface for $K_s \# (-K_s)$.  Tying $J$
  and $J'$, the Seifert surface of $K_s \# (-K_s)$ becomes a Seifert
  surface $E$ of genus $2r$ for~$L$.  $E$~consists of a single
  0-handle and $4r$ 1-handles $H_1,\ldots,H_{4r}$, where $H_i$ is the
  image of the~$H_{4r-i+1}$ under an obvious reflection, for $2\le i
  \le 2r$.  Joining the endpoints of the core of $H_i$ to their image
  under the reflection using disjoint arcs on the 0-handle of $E$ for
  $2\le i \le 2r$, we obtain $(2r-1)$ disjoint circles
  $\alpha_2,\ldots, \alpha_{2r}$ on~$E$.  See
  Figure~\ref{figure:infected-connected-sum}.

  \begin{figure}[ht]
    \begin{center}
      \includegraphics{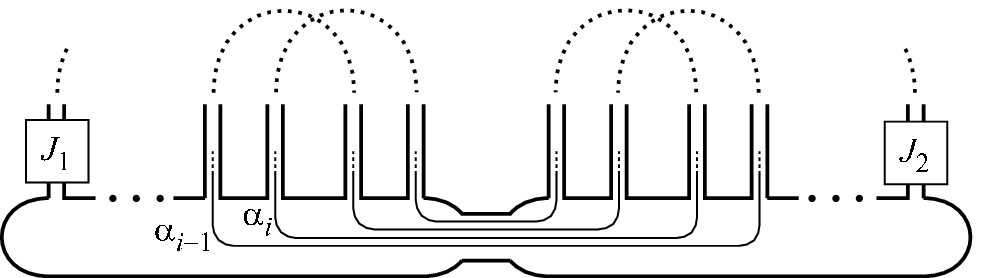}
    \end{center}
    \caption{}
    \label{figure:infected-connected-sum}
  \end{figure}
  
  The union of the $\alpha_i$ is a smoothly slice link, being the
  connected sum of a link and its mirror image.  Thus there are
  disjoint 2-disks $D_1, \ldots, D_{2r-1}$ smoothly embedded in $D^4$
  such that $\partial D_i = \alpha_i$.  Since the Seifert form defined
  on $E$ vanishes at $(\alpha_i, \alpha_j)$, one can do ambient
  surgery on $E$ along the $\alpha_i$, using the disks $D_i$ in $D^4$,
  as in~\cite{Levine:1969-1}.  This produces a genus one surface in
  $D^4$ with boundary~$L$.  Therefore $g_*^s(L) \le 1$.  This
  completes the proof of~(P2).

  Since $L$ shares a Seifert matrix with $K_s \# (-K_s)$ which is a
  ribbon knot, $L$ is algebraically slice.  From this (P3) follows.

  It is easily seen that $\Delta_K(t)=\Phi_n(t)^{2g}$.  Since $n$ has
  been chosen to be divisible by three distinct primes, (P4) holds due
  to a result of Livingston~\cite{Livingston:2002-1}.
\end{proof}

Let $\mathcal{C}$ be the smooth knot concordance group.

\begin{lemma}
  \label{lemma:main-example}
  Suppose $\{\Phi_\alpha\colon \mathcal{C} \to \Z\}$ is a finite
  collection of group homomorphisms satisfying $|\Phi_\alpha(-)| \le
  f_\alpha(g^s_*(-))$ for some real-valued function~$f_\alpha$.  Then,
  there are knots $J$ and $J'$ such that our $K$ satisfies the
  followings:
  \begin{enumerate}
  \item $g_*^h(K) = g_*^t(K)=g_*^s(K)=g$.
  \item $\Phi_\alpha(K)=0$ for each $\Phi_\alpha$.
  \end{enumerate}
\end{lemma}

Note that the Ozsv\'ath-Szab\'o $\tau$-invariant
\cite{Ozsvath-Szabo:2003-1} and Rasmussen $s$-invariant
\cite{Rasmussen:2004-1} can be viewed as homomorphisms of
$\mathcal{C}$ giving lower bounds of~$g_*^s$.  Therefore, from
Lemma~\ref{lemma:main-example}~(2), it follows that $J$ and $J'$ can
be chosen in such a way that $K$ has vanishing $\tau$- and
$s$-invariants.

\begin{proof}[Proof of Lemma~\ref{lemma:main-example}]
  Let $K'$ be the connected sum of $g$ copies of $K_s \# (-K_s)$.  By
  Cheeger-Gromov~\cite{Cheeger-Gromov:1985-1}, there is a universal
  bound $C$ of the $\rho$-invariants of the zero-surgery manifold
  $M_{K'}$ of $K'$, i.e., $|\rho(M_{K'},\phi')| \le C$ for any
  homomorphism $\phi'$ of~$\pi_1(M_{K'})$.

  Following \cite{Cochran-Orr-Teichner:2002-1}, for a knot $J$, let
  \[
  \rho(J)=\int_{S^1} \sigma_J(\omega) \, d\omega
  \]
  be the integral of the knot signature function
  \[
  \sigma_J(\omega)=\sign((1-\omega)S+(1-\bar \omega)S^T)
  \]
  over the unit circle $S^1$ normalized to unit length, where $S$ is a
  Seifert matrix of~$J$.

  We claim that there are two knots $J$ and $J'$ such that
  {
    \begin{enumerate}
    \item[(i)] $|\rho(J)| \ge C+4g$,
    \item[(ii)] $|\rho(J')| \ge C+4g+g\cdot |\rho(J)|$, and
    \item[(iii)] $\Phi_\alpha(K_s(\eta,J)) =
      \Phi_\alpha(K_s(\eta,J'))$ for each~$\Phi_\alpha$.
    \end{enumerate}
  }
  To prove the claim, we consider a sequence $\{J_i\}$ of knots
  constructed inductively as follows.  Let $J_0$ be a knot with
  $|\rho(J_0)| \ge C+4g$.  Assuming $J_{i}$ has been chosen, let
  $J_{i+1}$ be a knot satisfying
  \[
  |\rho(J_{i+1})| \ge C+4g+g\cdot |\rho(J_i)|.
  \]
  For example, one can choose as $J_i$ the connected sum of
  sufficiently many copies of any knot with nonvanishing $\rho$, e.g.,
  the trefoil knot, since $\rho$ is additive under connected sum.

  Since
  \[
  g_*^s(K_s(\eta,J_i)) \le g(K_s(\eta,J_i)) \le g(K_s)
  \]
  where $g(-)$ denotes the 3-genus (Seifert genus), there is an upper
  bound, say $M_\alpha$, of $f_\alpha\big(g_*^s(K_s(\eta,J_i))\big)$,
  i.e., $f_\alpha\big(g_*^s(K_s(\eta,J_i))\big) \le M_\alpha$ for any
  $J_i$.  Since
  \[
  |\Phi_\alpha(K_s(\eta,J_i))| \le f_\alpha\big(g_*^s(K_s(\eta,J_i))\big)
  \]
  by our hypothesis, it follows that $|\Phi_\alpha(K_s(\eta,J_i))|$ is
  bounded by $M_\alpha$.  Therefore the function $\Z\to
  \Z^{|\{\Phi_\alpha\}|}$ given by
  \[
  i \to (\Phi_\alpha(K_s(\eta,J_i)))_\alpha
  \]
  has finite image.  It follows that for some $i<j$,
  $\Phi_\alpha(K_s(\eta,J_i)) = \Phi_\alpha(K_s(\eta,J_j))$ for
  each~$\Phi_\alpha$.  Choosing $J=J_i$ and $J'=J_j$, the claim
  follows.  (Indeed our argument shows that there are infinitely many
  pairs $(J,J')$ satisfying the desired properties.)

  Recall that our $K$ is given by
  \[
  K = \mathop{\#}^g  \big(K_s(\eta, J) \# -(K_s(\eta,J')) \big).
  \]
  By (iii) above, $\Phi_\alpha$ vanishes at $K_s(\eta, J) \#
  -(K_s(\eta,J'))$.  It follows that $\Phi_\alpha(K)=0$ for
  each~$\Phi_\alpha$.  This proves the second conclusion of the lemma.

  To prove the first conclusion, it suffices to show that $g_*^h(K)
  \ge g$ by the property (P2) above.  Suppose $g_*^h(K) < g$.  By
  Proposition~\ref{proposition:bordism-from-slice-surface}, there is a
  4-manifold $V$ bounded by $M_K$ such that $\beta_2(V) < 2g$ and
  $\phi\colon\pi_1(M_K) \to \Z$ factors through~$\pi_1(V)$.

  Letting $\Gamma=\Z$, $\sR=\Q[t,t^{-1}]$, and $\K=\Q(t)$, we will
  apply
  Theorem~\ref{theorem:lower-bound-of-minimal-second-betti-number}~(2)
  to obtain a new coefficient system~$\phi_1$ which is a lift of
  $\phi$.  The conditions required in
  Theorem~\ref{theorem:lower-bound-of-minimal-second-betti-number}~(2)
  are verified as follows.  It is well-known that $A_K = H_1(M_K;\sR)$
  is always $\sR$-torsion.  We claim that $A_K$ is not generated by
  $\beta_2(V)$ elements.  Since the Alexander module is additive under
  connected sum and the knots $K_s(\eta, J)$ and $K_s(\eta,J')$ share
  the Alexander module with $K_s$, we have $A_K = \bigoplus^{2g}
  A_{K_s}$.  Since $A_{K_s}$ is nontrivial and $\beta_2(V)<2g$, $A_K$
  is never generated by $\beta_2(V)$ elements as claimed, by appealing
  to the structure theorem of finitely generated modules over
  $\Q[t,t^{-1}]$.

  Therefore, by applying
  Theorem~\ref{theorem:lower-bound-of-minimal-second-betti-number}~(2)
  and then~(1), it follows that there is a nontrivial homomorphism
  $h\colon A_K \to \K/\sR$ that gives rise to a homomorphism
  \[
  \phi_1\colon \pi_1(M_K) \to (\K/\sR) \rtimes \Gamma
  \]
  such that 
  \[
  |\rho(M_K,\phi_1)| \le 2\beta_2(V) < 4g. \tag*{(*)}
  \]
  
  Note that $K$ can be viewed as a knot obtained from $K'$ by tying
  $J$ and $-J'$ $g$ times.  So, from \cite[Proposition
  3.2]{Cochran-Orr-Teichner:2002-1} it follows that for some
  $\phi'\colon \pi_1(M_{K'}) \to (\K/\sR)\rtimes \Gamma$,
  \[
  \rho(M_K,\phi_1) = \rho(M_{K'},\phi') + \sum_{i=1}^g n_{i} \rho(J) -
  \sum_{i=1}^g m_i \rho(J').
  \]
  Here $n_{i}=0$ if the $(2i-1)$-st factor of $A_K = \bigoplus^{2g}
  A_{K_s}$ is contained in the kernel of $h$, and $n_i=1$ otherwise.
  The $m_i$ are determined similarly by the behaviour of the $(2i)$-th
  factor of~$A_K$.

  Since $h$ is a nontrivial homomorphism of $A_K$, at least one $n_i$
  or $m_i$ is nonzero.  If $m_i=0$ for all $i$, then since $n_i\ne 0$
  for some $i$, we have
  \[
  |\rho(M_K,\phi_1)| \ge |\rho(J)| - |\rho(M_{K'}, \phi')| \ge
  (4g+C)-C = 4g
  \]
  by (i) above.  It contradicts~$(*)$.  Therefore $m_i\ne 0$ for
  some~$i$.  In this case, by (ii) above, we have
  \[
  |\rho(M_K,\phi_1)| \ge |\rho(J')| - g\cdot |\rho(J)|-|\rho(M_{K'},
  \phi')| \ge (4g+C)-C = 4g.
  \]
  It again contradicts~$(*)$.  This shows that $g_*^h(K) \ge g$.
\end{proof}

\begin{remark}
  It can be easily seen that our construction produces infinitely many
  knot types of $K$.  In fact, one can use infinitely many knot types
  as our~$K_s$, $J$, and~$J'$.
\end{remark}

\begin{remark}
  Using results in \cite{Cochran-Orr-Teichner:1999-1}, it can be shown
  that the nonvanishing of the $\rho$-invariants we considered in the
  proof of Lemma~\ref{lemma:main-example} implies that our $K$ is not
  topologically slice.  Our result generalizes this.  In fact, our
  construction can be used to construct $K$ which is $(1)$-solvable
  but not $(1.5)$-solvable, in the sense of
  \cite{Cochran-Orr-Teichner:1999-1}.  It would be an interesting
  question whether there are $(h)$-solvable knots with topological
  slice genus $g$ for any $h\in \frac12 \Z$ and any $g>1$.
\end{remark}

We finish this section with an observation on the failure of an
attempt to extract information on the minimal genus for our example
using previously known results.  In \cite{Kervaire-Milnor:1961-1,
  Hsiang-Szczarba:1971-1, Rokhlin:1971-1, Lee-Wilczynski:1997-1,
  Lee-Wilczynski:2000-1} lower bounds of the topological minimal genus
are obtained for a homology class $\sigma\in H_2(X)$ in a topological
4-manifold $X$ which is closed or has boundary consisting of homology
sphere components.  When $X$ is simply connected,
\cite{Kervaire-Milnor:1961-1}~provides an obstruction to being
represented by a locally flat sphere, i.e., minimal genus $\ge 1$,
based on the Rokhlin theorem.  When $H_1(X)=0$,
\cite{Hsiang-Szczarba:1971-1, Rokhlin:1971-1, Lee-Wilczynski:1997-1,
  Lee-Wilczynski:2000-1} provides higher lower bounds of the following
form:
\[
2 \cdot(\text{minimal genus}) \ge -\beta_2(X)+ \max_{0 \le j <d}
\bigg| \operatorname{sign}(X)-\frac{2j(d-j)}{d^2} (\sigma \cdot
\sigma) \bigg|
\]
where $d$ is a positive integer such that $\sigma$ is contained in the
subgroup $d\cdot H_2(X)$.  (A more refined result of Lee-Wilczy\'nski
\cite[Theorem~2.1]{Lee-Wilczynski:2000-1} may potentially give further
lower bounds, however, computation seems infeasible when $H_1(X)\ne
0$.)

For an arbitrary 4-manifold $W$ with boundary and $\sigma\in H_2(W)$,
if $W$ embeds into a 4-manifold $X$ such that the above inequality
gives a lower bound for the image of $\sigma$ in $X$, then the result
is also a lower bound for $\sigma$ in~$W$.  However, when the
self-intersection of $\sigma$ is trivial, the above inequality gives
no information since $\beta_2(X) \ge |\operatorname{sign}(X)|$ for
any~$X$.  In particular, in the 4-manifold $W$ obtained by attaching a
2-handle to the 4-ball along the zero-framing of the knot $K$
constructed in this section, the generator $\sigma\in H_2(W)$ has
vanishing self-intersection so that the minimal genus is not detected
in this way.  Our results show that the minimal genus for $\sigma$ in
$W$ is exactly~$g$.

\bibliographystyle{amsalphaabbrv}
\renewcommand{\MR}[1]{}

\bibliography{research}

\end{document}